\documentclass[12TP,draft]{article}

\begin{document}

\begin{center}
\LARGE\noindent\textbf{On Hamiltonian and Hamilton-connected digraphs}\\

\end{center}

\begin{center}
\noindent\textbf{S.Kh.                                                                                                                                                                                                                                                                                                                                                                                                     Darbinyan }\\

Institute for Informatics and Automation Problems, Armenian National Academy of Sciences,

 P. Sevak 1, Yerevan 0014, Armenia

e-mail: samdarbin @ ipia.sci.am\\

A translation from Russian of the paper of Darbinyan ("O gamiltonovix i silno gamiltono-svyznix orgrafax", Akad. Nauk Armyan SSR Dokl. 91(1) 1990, 3-6). 
\end{center}
\begin{center}
\noindent\textbf{Abstract}\\
\end{center}

C. Thomassen in  \cite{[11]} suggested (see also \cite{[2]},  J. C.Bermond, C. Thomassen, Cycles in Digraphs - A survey, J. Graph Theory 5 (1981) 1-43, Conjectures 1.6.7 and 1.6.8) the following conjectures :

1. Every 3-strongly connected  digraph of order $n$ and with minimum degree at least $n+1$ is strongly  Hamiltonian-connected.

2. Let $D$ be a 4-strongly connected digraph of order $n$ such that the sum of the degrees of any pair of non-adjacent vertices is at least $2n+1$. Then $D$ is  strongly  Hamiltonian-connected.

We disprove   Conjecture 1  and prove two results which provide some support for Conjecture 2. The main goal of this article is to present the detailed proofs of these results (in English).\\


\noindent\textbf{1. Introduction}\\

 A cycle (path) of a directed graph (digraph) $D$ is called Hamiltonian if it  includes all the vertices of $D$. A digraph $D$ is Hamiltonian  if it contains a Hamiltonian cycle.

 Let us recall the following four well-known degree conditions (Theorems 1-4) for existence of a Hamiltonian cycles in digraph.  \\

 \textbf{Theorem 1.1} (Nash-Williams \cite{[9]}). {\it Let $D$ be a digraph of order $n\geq 2$ such that for every vertex $x$, $d^+(x)\geq n/2$ and $d^-(x)\geq n/2$, then $D$ is Hamiltonian.}\\

 \textbf{Theorem 1.2} (Ghouila-Houri \cite{[5]}). {\it Let $D$ be a strongly connected digraph of order $n\geq 2$. If $d(x)\geq n$ for all vertices $x\in V(D)$, then $D$ is Hamiltonian.}\\

 \textbf{Theorem 1.3} (Woodall \cite{[12]}). {\it Let $D$ be a digraph of order $n\geq 2$. If $d^+(x)+d^-(y)\geq n$ for all pairs of vertices $x$ and $y$ such that there is no arc from $x$ to $y$, then $D$ is Hamiltonian.}\\
 
\textbf{Theorem 1.4} (Meyniel \cite{[8]}). {\it  Let $D$ be a strongly connected digraph of order $n\geq 2$. If $d(x)+d(y)\geq 2n-1$ for all pairs of non-adjacent vertices in $D$, then $D$ is Hamiltonian.}\\

The existence of a Hamiltonian path with prescribed ends, is one of the extensions (generalizations) of the Hamiltonian cycle problem. A digraph $D$ is  strongly Hamiltonian-connected (respectively, weakly Hamiltonian-connected) if for any pair of distinct vertices $x$ and $y$  of $D$, there is a Hamiltonian path from $x$ to $y$ and  a Hamiltonian path  from $y$ to $x$ (respectively, there is a Hamiltonian path from $x$ to $y$ or  a Hamiltonian path  from $y$ to $x$).\\ 

As shown by Ghouila-Houri  \cite{[6]} for strongly Hamiltonian-connectedness  the analogous generalization of his theorem and, thus, of Meyniel's theorem is not true. 
He \cite{[6]} proved that every 2-strongly connected digraph of order $n\geq 2$  and with minimum degree at least $n+1$ is weakly Hamiltonian-connected. 
Furthermore, Meyniel's theorem cannot immediately be generalized, as there exist infinitely many 2-strongly connected tournaments that are not weakly Hamiltonian-connected \cite{[11]}.
A complete characterization of weakly Hamiltonian-connected tournaments was given by Thomassen \cite{[11]}.
 
Overbeck-Larisch \cite{[10]} considered the digraphs with condition of the type of  the condition Woodall's theorem, and the digraphs with condition of the type of  the condition Meyniel's theorem. She proved the following two theorems below.\\

\textbf{Theorem 1.5} (Overbeck-Larisch \cite{[10]}). {\it  Let $D$ be a digraph of order $n\geq 2$. Suppose that $d^+(x)+d^-(y)\geq n+1$ for each   pair of distinct vertices $x$ and $y$ such that there is no arc from $x$ to $y$, then $D$ is strongly Hamiltonian-connected.}\\

\textbf{Theorem 1.6} (Overbeck-Larisch \cite{[10]}). {\it Let $D$ be a 2-strongly connected digraph of order $n\geq 2$ with minimum degree at least $n+1$. Then $D$ is weakly Hamiltonian-connected.}\\

   A complete characterization of weakly Hamiltonian-connected tournaments was given by Thomassen \cite{[11]}. In \cite{[11]}, also was proved that every 4-strongly connected semicomplete digraph is strongly Hamilto-nian-connected and gave an infinite family of 3-strongly connected  tournaments with two vertices $x, y$, for which there is no   Hamiltonian path from $x$ to $y$. Thomassen also \cite{[11]} shows that for each $k\geq 1$ there exists a 2-strongly connected  non-strongly Hamiltonian-connected  digraph $D$ with minimum degree at least $|V(D)|+k$. In \cite{[11]}, the following two conjectures are given:\\

\noindent\textbf{Conjecture 1.7} (Thomassen \cite{[11]},  Conjecture 1.6.7 of \cite{[2]}). {\it Every 3-strongly connected digraph of order $n$ and with minimum degree at least $n+1$ is strongly Hamiltonian-connected.}\\

\noindent\textbf{Conjecture 1.8} (Thomassen \cite{[11]},  Conjecture 1.6.8 of \cite{[2]}). {\it Let $D$ be a  4-strongly connected  digraph  of order $n$ such that the sum of the degrees of any pair of non-adjacent vertices at least $2n+1$. Then $D$ is strongly Hamiltonian-connected.}\\

In this paper we disprove Conjecture 1.7 and prove two results which provide some support for  Conjecture 1.8. Below we will give the detailed proofs of these  results.\\

\noindent\textbf{2. Notation and Terminology }\\

We shall assume that the reader is familiar with the standard  terminology on graphs and directed graphs (digraphs) and refer  to  \cite{[1]} for terminology not discussed here.

In this paper we shall consider finite digraphs without loops and multiple arcs. For a digraph $D$, we denote by $V(D)$ the vertex set of $D$ and by  $A(D)$ the set of arcs in $D$. The arc of a digraph $D$ directed from $x$ to $y$ is denoted by  $xy$. For a pair of subsets $A$ and  $B$ of $V(D)$  we define $A(A\rightarrow B):= \{xy\in A(D) : x\in A, y\in B\}$,  $A(A,B):=A(A\rightarrow B)\cup A(B\rightarrow A)$. If $x\in V(D)$ and $A=\{x\}$, we often write $x$ instead of $\{x\}$. 

The out-neighbourhood of a vertex $x\in V(D)$ is the set $N^+_D(x):=\{y\in V(D) : xy\in A(D)\}$ and $N^-_D(x):=\{y\in V(D) : yx\in A(D)\}$ is the in-neighbourhood of $x$. Similarly, if $A\subseteq V(D)$ then $N^+_D(x,A):=\{y\in A : xy\in A(D)\}$ and $N^-_D(x,A):=\{y\in A : yx\in A(D)\}$. 
The out-degree of $x$ is $d^+_D(x):=|N^+_D(x)|$ and $d^-_D(x):=|N^-_D(x)|$ is the in-degree of $x$. Similarly, $d^+_D(x,A):=|N^+_D(x,A)|$ and $d^-_D(x,A):=|N^-_D(x,A)|$. The degree of a vertex $x$ in $D$ is $d_D(x):=d^+_D(x)+d^-_D(x)$. We usually drop the subscript $D$ if this is unambiguous. 

The subdigraph of $D$ induced by a subset $A$ of $V(D)$ is denoted by $D\langle A\rangle$, or $\langle A\rangle$ for brevity. If $A\subset V(D)$, then we denote by $D-A$ the subdigraph $\langle V(D)\setminus A \rangle$. 

For integers $a$ and $b$, $a\leq b$, let $[a,b]$  denote the set of all integers which are not less than $a$ and are not greater than $b$. 

The path (respectively, the cycle) consisting of the distinct vertices $x_1,x_2,\ldots ,x_n$ ($n\geq 2 $) and the arcs $x_ix_{i+1}$, \,$i\in [1,n-1]$  (respectively, $x_ix_{i+1}$, \,$i\in [1,n-1]$ and $x_nx_1$), is denoted by $x_1x_2\ldots x_n$ (respectively, $x_1x_2\ldots x_nx_1$). 
The path $x_1x_2\ldots x_n$ is called an $(x_1,x_n)$-path or a path from $x_1$ to $x_n$. The cycle $x_1x_2\ldots x_nx_1$ (respectively, the path  $x_1x_2\ldots x_n$) of $D$ is Hamiltonian if $n=|V(D)|$.
 For a cycle  $C_k:=x_1x_2\ldots x_kx_1$, the indices are taken modulo $k$, i.e. $x_s=x_i$ for every $s$ and $i$ such that  $i\equiv s\, \hbox {(mod} \,k)$. If $C$ is a cycle containing vertices $x$ and $y$, $C[x,y]$ denotes the subpath of $C$ from $x$ to $y$.

 A digraph $D$ is strongly connected, or strong for brevity, (respectively, unilaterally connected) if for every pair  of distinct vertices $x$, $y$ of $D$ there exists an $(x,y)$-path and a $(y,x)$-path (respectively, an $(x,y)$-path or a $(y,x)$-path). 
  A digraph $D$ is $k$-strongly connected (or $k$-strong, $k\geq 1$) if $|V(D)|\geq k+1$ and $D-A$ is strong for any set $A$ of at most $k-1$ vertices. 
By Menger's theorem, this is equivalent to the property that for any ordered pair of distinct vertices $x,y$ there are $k$ internally disjoint paths from $x$ to $y$. 
 A strong component  of a digraph $D$ is a maximal induced strong subdigraph of $D$. Two distinct vertices $x$ and $y$ in $D$ are adjacent if $xy\in A(D)$ or $yx\in A(D) $ (or both). 

  The converse digraph of a digraph $D$ is the digraph obtained from $D$ by reversing the directions of all arcs of $D$.  
 \\

\noindent\textbf{3. Non-strongly Hamiltonian-connected 3-strong digraphs with large minimum degree }\\
 
In this section using the construction (as well as the converse construction) of M. Overbeck-Larisch \cite{[10]} we will disprove  Conjecture 1.7.

\noindent\textbf{The construction of M. Overbeck-Larisch.} Let $D$ be a digraph of order $n+1\geq 5$  and  $u$ and $v$ be arbitrary two distinct vertices of $D$. Construct (see M. Overbeck-Larisch \cite{[10]}) a new digraph $H_D(u,v)$ of order $n$ with
 $$V(H_D(u,v)):=(V(D)\setminus \{u,v \})\cup \{z_0\} \quad ( z_0 \quad \hbox {a new vertex})$$
 and
$$A(H_D(u,v)):=A(D-\{u,v\})\cup \{z_0y : y\in N^+_{D-\{v\}}(u)\}\cup \{yz_0 : y\in N^-_{D-\{u\}}(v)\}.$$
It is not difficult to see that $|V(H_D(u,v))|=n$ and  $d_{H_D(u,v)}(z_0)=d^{+}_{D-\{v\}}(u)+d^-_{D-\{u\}}(v)$.\\ 

\noindent\textbf{Lemma 3.1}. {\it Let $D$ be a digraph of order $n+1\geq 5$. If $D$ is $k$-strong ($k\geq 3$), then for any two distinct vertices $u$, $v$ of $D$, $H_D(u,v)$ is $(k-1)$-strong.}

\noindent\textbf{Proof of Lemma 3.1}. Assume that $D$ is $k$-strong but $H:=H_D(u,v)$ is not $(k-1)$-strong. 
Then there are $k-2$ distinct vertices $x_1,x_2,\ldots ,x_{k-2}\in V(H)$ such that  $ H-\{x_1,x_2,\ldots ,x_{k-2}\} $ is not strong, i.e., for some two distinct vertices $x, y \in V(H)\setminus \{x_1,x_2,\ldots ,x_{k-2}\}$ in $H- \{x_1,x_2,\ldots ,x_{k-2}\} $ there is no $(x,y)$-path. 
From $k$-strong connectedness of $D$ it follows that  $z_0\notin \{x_1,x_2,\ldots ,x_{k-2}\}$ (for otherwise, if $z_0=x_i$,
say $z_0=x_{k-2}$, then $D-\{x_1,x_2,\ldots , x_{k-3},u,v\}$ is not strong, which is a contradiction), and by Menger's theorem in $D$ there are $k$ internally disjoint $(x,y)$-paths. Assume that $x$ and $y \in V(H)\setminus \{z_0\}$. 
Then one of these paths necessarily has  the following form
  $xu_1\ldots u_juu_{j+1}\ldots u_ly$; and another of these paths necessarily has the following form $xv_1\ldots v_rvv_{r+1}\ldots v_qy$, and both   are in $D-\{x_1,x_2,\ldots ,x_{k-2}\}$. 
Therefore   $xv_1\ldots v_rz_0u_{j+1}\ldots u_ly$ is an $(x,y)$-path in  $H- \{x_1,x_2,\ldots ,x_{k-2}\}$, a contradiction. Now assume that  $x=z_0$ (the argument for $y=z_0$ is similar). 
Then, since $D$ is $k$-strong,  in $D-\{x_1,x_2,\ldots ,x_{k-2},v\}$ there is a path $ua_1a_2\ldots a_qy$, and therefore $z_0a_1a_2\ldots a_qy$ is an $(x,y)$-path in $ H-\{x_1,x_2,\ldots ,x_{k-2}\} $. 
This contradicts the assumption that in $ H-\{x_1,x_2,\ldots ,x_{k-2}\}$ there is no $(x,y)$-path. Therefore, $H$ is ($k-1$)-strong. Lemma 3.1 is proved. \fbox \\\\

\noindent\textbf{The converse construction of M. Overbeck-Larisch construction.} Let $H$ be a digraph of order $n\geq 4$ and let $z_0$ be an arbitrary vertex of $H$. Now we define a digraph $D_H(z_0)$ as follows: 
$$V(D_H(z_0):=(V(H)\setminus \{z_0\}) \cup\{u,v \} \quad (u,v \quad \hbox {are new vertices}),$$
$$A(D_H(z_0)):=A(H- \{z_0\}) \cup \{uv,vu\} \cup \{xu,vx :x\in V(H)\setminus \{z_0\}\}$$ $$\cup \{ux :z_0x\in A(H)\}\cup  \{xv :xz_0\in A(H)\}.$$ 
Note that $D_H(z_0)$ has $n+1$ vertices. \\

\noindent\textbf{Lemma 3.2}. {\it Let $H$ be a digraph of order $n\geq 4$. If $H$ is $k$-strong ($k\geq 2$), then for every vertex $z$ of $H$, $D_H(z)$ is $(k+1)$-strong.}

\noindent\textbf{Proof of Lemma 3.2}. Assume that the digraph $H$ is $k$-strong, but for some vertex $z_0\in V(H)$, $D:=D_H(z_0)$ is not $(k+1)$-strong. 
Then there are $k$ distinct vertices $x_1,x_2,\ldots ,x_{k}\in V(D)$ such that $D- \{x_1,x_2,\ldots ,x_{k}\}$ is not strong, i.e., for some two distinct vertices $x, y \in V(D)\setminus \{x_1,x_2,\ldots ,x_{k}\}$ in $D- \{x_1,x_2,\ldots ,x_{k}\}$ there is no $(x,y)$-path. Note that $\{u,v\}\not\subseteq \{x_1,x_2,\ldots ,x_{k}\}$, since $H$ is $k$-strong.

Let $\{x_1,x_2,\ldots ,x_{k}\} \subset V(D)\setminus \{u,v \}$. 
Then it is easy to see that $\{ x,y \}\not= \{u,v \}$ (for otherwise, $xy\in A(D)$), and $\{ x,y \}\not\subset V(D)\setminus \{u,v \}$ (for otherwise, $xuvy$ is an $(x,y)$-path in $D-\{x_1,x_2,\ldots ,x_k\}$). 
This mean that either  $x \in  \{u,v \}$ or $y\in \{u,v\}$. Now it is easy to see that in  both cases  in  $D- \{x_1,x_2,\ldots ,x_k\}$ there is an $(x,y)$-path.

Let now  $\{x_1,x_2,\ldots ,x_k\}\notin \subseteq V(D)\setminus \{u,v \}$. 
Without loss of generality we may assume that
 $\{x_1,x_2,\ldots ,$ $x_{k-1}\}\subset V(D)\setminus \{u,v \}$, i.e., $x_k\in \{u,v\}$. 
If $x,y\in V(D)\setminus \{u,v\}$, then $H$ is not $(k+1)$-strong and in $H-\{x_1,x_2,\ldots ,x_{k-1}\}$ there exists a path $xv_1\ldots v_iz_0v_{i+1}\ldots v_ly$ since $H$ is $k$-strong. 
From this it is not difficult to see that in  $D-\{x_1,x_2,\ldots ,x_k\}$, if $x_k=v$, then $xu v_{i+1}\ldots v_ly$ is an $(x,y)$-path, and if $x_k=u$, then $xv_1\ldots v_ivy$ is an $(x,y)$-path. 
We may therefore  assume that either $x\in \{u,v \}$ or $y\in \{u,v \}$. We will  consider only the case $x\in \{u,v\}$ (the argument for $y\in \{u,v\}$ is similar). Then  $x=u$ and  $x_k=v$. 
In $H-\{x_1,x_2,\ldots ,x_{k-1}\}$ there is a path  $z_0v_{1}\ldots v_ry$ since $H$ is $k$-strong. 
Therefore, $uv_{1}\ldots v_ry$ is an $(x,y)$-path in $D-\{x_1,x_2,\ldots ,x_k\}$. 
Thus in all possible cases in $D-\{x_1,x_2,\ldots ,x_k\}$ there is an $(x,y)$-path, which contradicts the  assumption that $D$ is not $(k+1)$-strong. This   completes the proof of the lemma. \fbox \\\\ 
 
\noindent\textbf{Theorem 3.3.}  {\it Every $k$-strong ($k\geq 2$) digraph of order $n\geq 3$ which has $n-1$ vertices of degree at least $n$ is Hamiltonian if and only if any $(k+1)$-strong  digraph of order $n+1$ with minimum degree at least
 $n+2$ is strongly Hamiltonian-connected.}

\noindent\textbf{Proof of Theorem 3.3}. Suppose that every $k$-strong digraph of order $n\geq 3$  which has $n-1$ vertices of degree at least $n$ is Hamiltonian. Let $D$ be a  $(k+1)$-strong  digraph of order $n+1$  with  minimum degree at least $ n+2$. Let $u$ and $v$ be two arbitrary  distinct vertices of $D$. Consider the digraph $H:=H_D(u,v)$. Recall that $H$ has $n$ vertices. 
It is not difficult to see that  $d_H(x)\geq n$  for all $x\in V(H)\setminus \{z_0\}$ (note that $d_H(z_0)=d^{+}_{D-\{v\}}(u)+d^-_{D-\{u\}}(v)$). By Lemma 3.1, $H$ is $k$-strong. Therefore, by the our supposition, $H$ contains a Hamiltonian cycle,  which in turn implies that $D$ has a Hamiltonian $(u,v)$-path.

 Now suppose that every $(k+1)$-strong digraph of order $n+1$  with minimum degree at least $n+2$ is strongly Hamiltonian-connected. 
Let $H$ be a $k$-strong digraph of order $n$  whose $n-1$ vertices have degrees at least $n$. Let the vertex $z_0$ has the minimum degree in $H$.  Now we consider the digraph $D:=D_H(z_0)$. Note that $|V(D)|=n+1$ and
$$
 d_D(x)\geq \left\{ \begin{array}{lc} n+2,\quad \hbox{if}\quad x\in V(D)\setminus\{u,v\},    \\ n+k+1, \quad \hbox{if}  \quad x\in \{u,v\}.  \\ \end{array} \right.
$$
On the other hand, by Lemma 3.2, $D$ is $(k+1)$-strong and hence, by the supposition, in $D$ there is a Hamiltonian  $(u,v)$-path., which in turn . Therefore $H$ is Hamiltonian. Theorem 3.3  is proved. \fbox \\\\

\noindent\textbf{Theorem 3.4.}  {\it For every integer $n\geq 8$ there is a 2-strong non-Hamiltonian digraph  of order $n$  which has $n-1$  vertices of degrees  at least $n$.}

\noindent\textbf{Proof of Theorem 3.4}. We define a digraph $D$ on $n\geq 8$ vertices as follows:

$V(D):=\{x_0,x_1, \ldots , x_{n-4},y_1,y_2,y_3\}$ and
$$A(D):=\{y_iy_j : i \not=j\} \cup   \{x_ix_{i+1} : i\in [0, n-5]\} \cup  \{y_ix_j : i\in [1,  3], j\in [1,n-6]  \}  \cup  \{ x_ix_j/ 1 \leq j< i\leq n-4 \}$$
$$
\cup   \{x_{n-4}y_{i}, x_{n-6}y_i : i\in [1, 3]\} \cup \{x_ix_{n-5} : i\in [1, n-7]\} \cup  \{x_0x_{n-5}, x_{n-5}x_{0},x_{n-4}x_{0},x_{n-6}x_{n-4}\}.
$$
Observe that $d(y_i)=n$ for all $i\in [1,3]$, $d(x_0)=4$,  $d(x_j)=n+1$ for all $j\in [1,n-7]$, $d(x_{n-4})=n+1$, $d(x_{n-5})=2n-8$ and $d(x_{n-6})=n+4$. Therefore, $n-1$ vertices of  $D$ have degrees  at least $n$.\\ 

Using the following simple 

{\it Proposition:  Let $H$ be a strong digraph and let $x$ be a new vertex not in $H$. 
Let 
$H'$ be a digraph obtained from  $H$ by adding a new vertex $x$ and adding an arc from $x$ to a vertex of $H$ and an arc from a vertex of $H$ to $x$, then $H'$ also is strong;} 

it is not difficult to check that for each $z\in V(D)$ the digraph $D-z$ is strong, i.e., $D$ is 2-strong.  

Now we prove that $D$ is not Hamiltonian. Suppose, to the contrary, that $D$ is Hamiltonian. Let $C$ be an arbitrary Hamiltonian cycle in $D$.  It is not difficult to see that the path $C[x_1,x_{n-5}]$ necessarily has the following type $x_1x_2\ldots x_ix_{n-5}$, where $i\in [1,n-6]$.
 Hence, from the construction of $D$ it follows that the cycle $C$ does  not contain the arc $x_{n-5}x_0$. Therefore $C$ contains the arc $x_{n-4}x_0$ and  either the path $x_{n-5}x_{n-4}x_{0}$ or the path $x_{n-6}x_{n-4}x_{0}$. Now again by the construction of $D$ we obtain that the cycle $C$ does not contain the vertices $y_1,y_2,y_3$. This contradicts that $C$ is a Hamiltonian cycle in $D$. Theorem 3.4 is proved. \fbox \\\\

\noindent\textbf{Remark 3.5}. {\it Let $D$ be the digraph of order $n\geq 8$ that is defined in the proof of Theorem 3.4. If $xy\notin A(D)$, then the digraph $H$ obtained from $D$ by adding the arc $xy$, is Hamiltonian.} \\

From Theorems 3.3 and 3.4 it follows the following theorem which disproves Conjecture 1.7.\\

\noindent\textbf{Theorem 3.6.} {\it  For every integer $n\geq 9$ there is a 3-strong   non strongly Hamiltonian-connected digraph of order $n$ with minimum degree at least  $n+1$.}

\noindent\textbf{Proof}. We define a digraph  $D$ as  in the proof of Theorem 3.4. Then the digraph $D_D(x_0)$ is 3-strong (Lemma 3.2) and has the minimum degree at least $|V(D_D(x_0)|+1$. By Theorems 3.3 and 3.4, in $D_D(x_0)$ there is no Hamiltonian $(u,v)$-path.  \fbox \\\\

\noindent\textbf {Note added in proof} (for section 3). Later on, a characterization of weakly Hamiltonian-connected tournaments due to Thomassen has been generalized to several classes of generalizations of tournaments. \\

(i) Bang-Jensen, Guo and Volkmann \cite{[13]} gave a complete characterization of weakly Hamiltonian-connected locally semicomplete  digraphs. 

(ii) Bang-Jensen, Gutin and Huang \cite{[14]} gave a characterization of weakly Hamiltonian-connected extended tournaments (ordinary multipartite tournaments). 

(iii) A characterization of weakly Hamiltonian-connected semicomplete bipartite digraphs without 2-cycles was obtained by Bang-Jensen and Manoussakis \cite{[15]}.\\

 In \cite{[4]}, the author proved the following theorem. \\

 \noindent\textbf{Theorem 3.7} (\cite{[4]}). {\it Let $D$ be a 2-strong digraph of order $n\geq 9$ with minimum degree at least $n-4$. If  $n-1$ vertices of $D$ have degrees at least $n$, then $D$ is  Hamiltonian.}\\

In \cite{[4]}, given only the outline of the proof of Theorem 3.7. The proof of Theorem 3.7 is rather lengthy and involves much cases analysis. We put as a question to find a sort proof of this result.\\

As noted above the digraph $D$ in Theorem 3.4 has the minimum degree equal to four. It is natural to pose the following.\\

\noindent\textbf{Problem 1}. {\it Let $D$ be an arbitrary 2-strong digraph of order $n\geq 8$. Suppose that  $n-1$ vertices of $D$ have degrees at least $n$ and a vertex $x$ has degree at least $n-m$, where $5\leq m\leq n-5$. Find the maximum volume of $m$ for which the digraph $D$ is  Hamiltonian.}\\

Observe that in the proof of Theorem 3.4 (respectively, Theorem 3.6) the vertex-connectivity of $D$ is equal to two
(respectively, the vertex-connectivity of $D_D(x_0)$ is equal to three). Now we can reformulate Conjecture 1.7 due to Thomassen in the following two forms.\\

\noindent\textbf{Problem 2}. {\it Let $D$ be a  digraph of order $n\geq 8$ in which  $n-1$ vertices  have degrees at least $n$. Does there exists a integer $k\geq 3$ such that if $D$ is $k$-strong, then $D$ is  Hamiltonian (Find the least volume of $k$ if it exists).}\\

\noindent\textbf{Problem 3}. {\it Let $D$ be a  digraph of order $n\geq 9$ with minimum degree at least  $n+1$. 
 Does there exists a integer $k\geq 4$ such that if $D$ is $k$-strong, then $D$ is  strongly Hamiltonian-connected (Find the least volume of $k$ if it exists).}\\

\noindent\textbf{4. Some supports for Conjecture 1.8 }\\

 In this section we prove two results which provide some supports for Conjecture 1.8. 

We say that a digraph $D$ of order $n$ satisfies  condition $(M)$:  
If the sum of the degrees of any pair of non-adjacent vertices $x,y \in V(D)\setminus \{z_0\}$  at least $2n-1$, where $z_0$ is some vertex of $D$; and satisfies  condition $(N)$: If the sum of the degrees of any pair of non-adjacent vertices $x,y$ 
of $D$    at least $2n+1$. \\

\noindent\textbf{Theorem 4.1.}  {\it Any $k$-strong ($k\geq 1$) digraph of order $n\geq 8 $  satisfying  condition $(M)$ is Hamiltonian if and only if any $(k+1)$-strong  digraph of order $n+1 $ satisfying  condition $(N)$ is strongly Hamiltonian-connected.} 

\noindent\textbf{Proof of Theorem 4.1}. Suppose first that any $k$-strong digraph of order $n\geq 8$  satisfying  condition $(M)$ is Hamiltonian.  Let $D$ be a  $(k+1)$-strong  digraph of order $n+1\geq 9 $  satisfying  condition $(N)$, and let $u$ and $v$ be any two distinct vertices of $D$. We want to prove that in $D$ there is a Hamiltonian $(u,v)$-path. Consider the digraph $H:=H_D(u,v)$. It is not difficult to see that $|V(H)|=n$ and $d_H(x)+d_H(y)\geq d_D(x)+ d_D(y)-4 \geq 2n-1$  for any two non-adjacent distinct vertices $x,y\in V(H)- \{z_0\}$, i.e., $H$ satisfies  condition $(M)$.
 From Lemma 3.1 it follows that $H$ is k-strong. By the supposition, $H$ is Hamiltonian, which in turn implies that  in $D$ there is a Hamiltonian $(u,v)$-path.

Now suppose that any $(k+1)$-strong digraph of order at least 9 satisfying condition $(N)$ is strongly Hamiltonian-connected.
Let $H$ be a  $k$-strong  digraph of order $n\geq 8 $  satisfying  condition $(M)$. Now we consider the digraph  $D:=D_H(z_0)$. By the construction of $D$, it is easy to see that $D$ has $n+1$ vertices and $d_D(x)=d_H(x)+2$  for all $x\in V(H)-\{z_0 \}$. From this and the construction of $D$ it follows that if the 
  vertices $x, y$ of $D$ are non-adjacent, then  $x$ , $y\in V(D)\setminus \{u,v\}$ and
$$
d_D(x)+d_D(y)\geq d_H(x)+2+d_H(y)+2\geq 2(n+1)+1,
$$
i.e., $D$ satisfies  condition $(N)$. On the other hand, by Lemma 3.2, $D$ is $(k+1)$-strong. Thus we have that $D$ is 
$(k+1)$-strong and satisfies  condition $(N)$.
By the supposition, in $D$ there is a Hamiltonian $(u,v)$-path. Therefore, $H$ contains a Hamiltonian cycle. Theorem 4.1 is proved. \fbox \\\\

\noindent\textbf{Remark 4.2}. {\it To establish Conjecture 1.8 it suffices to show that any 3-strong digraph of order $n\geq 8$ 
satisfying condition (M) is Hamiltonian.}\\ 

The following well-known simple lemmas is the basis of the proof of Theorem 4.5 and other theorems on directed cycles and paths in digraphs.\\

\noindent\textbf{Lemma 4.3}  (\cite{[7]}). {\it Let $D$ be a digraph of order $n\geq 3$  containing a cycle $C_m$, $m\in [2,n-1]$ and  let $x\notin C_m$. If $d(x,C_m)\geq m+1$, then for every $k \in [2,m+1]$ $D$ contains a cycle $C_k$ including $x$.}  \\

\noindent\textbf{Lemma 4.4} (\cite{[3]}). {\it Let $D$ be a digraph of order $n\geq 3$ containing a path $P:=x_1x_2\ldots x_m$, where $m\in [2,n-1]$. Let $x$ be a vertex not contained in this path. 
If $d(x,P)\geq m+2$, then there is an $i\in [1,m-1]$ such that $x_ix,xx_{i+1}\in D$, i.e., $D$ contains a path $x_1x_2\ldots x_ixx_{i+1}\ldots x_m$ of length $m$ (we say that the vertex $x$ can be inserted into $P$ or the path $x_1x_2\ldots x_ixx_{i+1}\ldots x_m$ is an extended path obtained from $P$ with the vertex $x$).}  \\

In the proof of Theorem 4.5 we often will use the following definition:

\textbf{Definition}. {\it Let $P_0:=x_1x_2\ldots x_m$, $m\geq 2$, be an $(x_1,x_m)$-path in a digraph $D$ and let the vertices
$y_1,y_2,\ldots y_k$ are in $V(D)- V(P_0)$. For $i\in [1,k]$ we denote by $P_i$ an $(x_1,x_m)$-path in $D$ with vertex set $V(P_{i-1})\cup \{y_j\}$ (if it exists) such that the path $P_i$ is an extended path obtained from the path $P_{i-1}$ with some vertex $y_j$, 
where $y_j\notin V(P_{i-1}$). If $e+1$ is the maximum possible number of these paths $P_0, P_1,\ldots , P_e$, $e\in [0,k]$, then we say that $P_e$ is an extended path obtained from $P_0$ with vertices $y_1,y_2,\ldots , y_k$ as much as possible. Notice that $P_i$ for all $i\in [0,e]$ is an $(x_1,x_m)$-path of length $m+i-1$.}\\

\noindent\textbf{Theorem 4.5.}  {\it Let $D$ be a strong digraph of order $n\geq 3$. If $d(x)+d(y)\geq 2n-1$ for any two non-adjacent vertices $x,y\in V(D)-\{z_0\}$, where $z_0$ is  some vertex of $D$ (i.e., $D$ satisfies  condition $(M)$), then $D$ is Hamiltonian or contains a cycle of length $n-1$.}

\noindent\textbf{Proof of Theorem 4.5}. Suppose that Theorem 4.5 is not true. Let $C_m:=x_1x_2\ldots x_mx_1$ be a cycle  of maximum length in $D$. Put $B:=V(D)\setminus V(C_m)$. It is clear that $|B|=n-m\geq 2$ and $n\geq 5$. By the maximality of $C_m$ and Lemma 4.3, for each vertex $x\in B$ we have
$$d(x, V(C_m))\leq m . \eqno {(1)} $$ 

First it is convenient to prove the following two claims below.

 \noindent\textbf{Claim 1}. {\it The subdigraph $\langle B-\{z_0\} \rangle$  is unilaterally connected.}

 \noindent\textbf{Proof of Claim 1}. If any two distinct vertices of $B\setminus \{z_0\}$ are adjacent, then Claim 1 is true. Assume that  some distinct vertices $x$ and $y$ of $B\setminus \{z_0\}$ are not adjacent. Then, by (1) and the hypotheses of the theorem, we have
$$
 2n-1\leq d(x)+d(y)\leq \left\{ \begin{array}{lc} 2m+d(x,B\setminus \{z_0\})+d(y,B\setminus \{z_0\})+d(z_0,\{x,y\}),\quad \hbox{if}\quad  z_0\in B ,    \\ 2m+d(x,B)+d(y,B),  \quad \hbox{otherwise } . \\ \end{array} \right. 
$$
Hence
$$
d(x,B\setminus \{z_0\})+d(y,B\setminus \{z_0\})\geq 2(n-m-1)-3, \, \hbox{if} \, z_0\in B,
$$
$$
d(x,B)+d(y,B)\geq 2(n-m)-1, \, \hbox{if}\, z_0\notin B.
$$
Now it is easy to see that in $\langle B\setminus \{z_o\} \rangle$ there is a path of length two with end-vertices $x$ and $y$. Claim 1 is proved.  \fbox \\\\

 \noindent\textbf{Claim 2}. {\it At least two distinct vertices of $C_m$ are adjacent to some vertices of $B\setminus \{z_o\}$.} 

\noindent\textbf{Proof of Claim 2}. Suppose that Claim 2 is not true. Without loss of generality, we may assume that
$$
A(\{ x_2,x_3,\ldots , x_m\}, B\setminus \{z_o\})=\emptyset. \eqno {(2)}
$$
From this and the  hypotheses of the theorem it follows that 
$$
d(x_i)+d(y)\geq 2n-1 \eqno (3)
$$
for every pair of vertices $x_i\in \{ x_2,x_3,\ldots , x_m\}\setminus \{z_0\}$ and $y\in B\setminus \{z_o\}$.

We distinguish two cases, according as $z_0$ is in $B$ or not.

\noindent\textbf{Case 1}. $z_0\in B$.
  
Assume first that the subdigraph $\langle V(C_m) \rangle$ is not a complete digraph. Then there is a vertex $x_k$, $k\in [2,m]$, such that $d(x_k,C_m)\leq 2m-3$. Now using (2) and (3), for $x_k$ and for every $y\in B\setminus \{z_0\}$ we obtain
$$
2n-1\leq d(x_k)+d(y) \leq 2m-3+d(x_k,\{z_0\})+2(n-m-2)+d(y,\{x_1,z_0\})=$$ 
$$2n-7+d(x_k,\{z_0\})+d(y,\{x_1,z_0\}).
$$
Therefore 
$$
d(x_k,\{z_0\})+d(y,\{x_1,z_0\})\geq 6 \quad \hbox{and} \quad d(x_k,\{z_0\})=d(y,\{z_0\})=d(y,\{x_1\})=2. \eqno (4)
$$
 Similar to (4), we will  obtain that every vertex of $\{x_2,x_3,\ldots, x_m\}$ is adjacent to $z_0$.  From the maximality of $C_m$ it follows that $z_0x_2\notin A(D)$ because of (4). Then $x_2z_0\in A(D)$, $k\geq 3$ and there is a $j\in [2,k-1]$ such that $x_jz_0$, $z_0x_{j+1}\in A(D)$, which is a contradiction, since $x_1x_2 \ldots x_jz_0x_{j+1} \ldots x_mx_1$ is a cycle of length $m+1$.

Now assume that $\langle V(C_m) \rangle$ is  complete. Then $D$ also contains the cycle $x_1x_mx_{m-1}\ldots x_2x_1$. Using (2) and (3), by the same arguments as  above, it is not difficult to show that 
$$
d(y,\{x_1,z_0\})+d(x_2,\{z_o\})\geq 5.
$$ 
Therefore, $D$ contains the path $x_1yz_0x_2$ or the path $x_2z_0yx_1$, where $y\in B\setminus \{z_0\}$. In each case in $D$ there is a cycle of length $m+2$, which is a contradiction.\\

\noindent\textbf{Case 2}. $z_0\notin B$, i.e., $z_0\in V(C_m)$.

Then using (2), we obtain  that $d(y)+d(x_i)\leq 2n-2$ for every pair of vertices $x_i$ ($i\in [2,m]$) and $y\in B$. This together with (3) implies that $x_i=z_0$, i.e., $m=2$. Since $D$ is strong and $\langle B \rangle$ is unilaterally connected (Claim  1), it is not difficult to see that $\langle B \rangle$ is strong and precisely one vertex of $B$ (say $x$) is adjacent to $x_1$ (for otherwise $D$ contains a cycle of length at least three). 
Therefore, $d(x_1,B\setminus \{x\})=0$ and for $x_1$ and  $u\in B\setminus \{x\}$ we have $d(x_1)+d(u)\leq 4+2(n-m-1)=2n-2$, which is a contradiction. This contradiction completes the proof of Claim 2.  \fbox \\\\

 Let $D_1, D_2,\ldots  , D_s$ ($s\geq 1$) be the strong components of the subdigraph $\langle B \rangle$ (or of the subdigraph $\langle B \setminus \{z_0\} \rangle$ if $ \langle B \rangle$ is not unilaterally connected) labeled in such a way that (by Claim 1 it is possible)
$$
A(V(D_i) \rightarrow V(D_{i+1})) \not= \emptyset \quad \hbox {and} \quad A(V(D_j) \rightarrow V(D_{i}))=\emptyset, \quad 1 \leq i < j\leq s. \eqno (5)
$$

 From Claim 1 it follows that if  $ \langle B \rangle$ is not unilaterally connected, then  $z_0\in B$, $$
A(z_0 \rightarrow V(D_{1})) = A(V(D_s) \rightarrow z_0)=\emptyset \eqno (6)
$$
and  the subdigraphs $ \langle V(D_i)\cup \{z_0\} \rangle$ are not strong for all $i\in [1,s]$.
 
Since $D$ is strong, from (5) and (6) it follows that
$$ 
A(V(C_m) \rightarrow V(D_{1})) \not= \emptyset \quad \hbox {and} \quad A(V(D_s) \rightarrow V(C_{m}))\not= \emptyset. \eqno (7)
$$
This together with Claim 2 imply that there are two distinct vertices $x_r, x_l\in V(C_m)$ such that either
$$ 
(i) \quad A(V(D_s) \rightarrow x_r)\not=\emptyset \quad \hbox {and} \quad A(x_l \rightarrow \cup_{i=1}^s V(D_{i}))\not=\emptyset \quad \hbox{or}
$$
$$ (ii) \quad A(x_l  \rightarrow V(D_1))\not=\emptyset \quad \hbox {and} \quad A( \cup_{i=1}^s V(D_{i})\rightarrow x_r)\not=\emptyset . 
$$
Assume that (i) holds (for the case (ii) we can apply the same arguments). By Claim 1 and the maximality of the cycle $C_m$, we have $|C_m[x_l,x_r]|\geq 3$ and $A(x_{r-1} \rightarrow \cup_{i=1}^s V(D_{i}))=\emptyset$. Now we can choose a vertex $x_a\in C[x_l,x_{r-2}]$ such that 
$$
 A(x_a \rightarrow \cup_{i=1}^s V(D_{i}))\not=\emptyset \quad \hbox {and} \quad A(C[x_{a+1},x_{r-1}] \rightarrow \cup_{i=1}^s V(D_{i}))=\emptyset. 
$$
Let $x_au\in D$, where $u\in D_k$ and $k$ is as small as possible.  
   By Claim 1 and the maximality of $C_m$ we have 
$$
A( \cup_{i=k}^s V(D_{i})\rightarrow x_{a+1})=\emptyset.
$$
 Therefore there is a vertex $x_b\in C[x_{a+2},x_r]$ such that 
$$
A( \cup_{i=k}^s V(D_{i})\rightarrow x_{b})\not=\emptyset \quad \hbox{and} \quad A( \cup_{i=k}^s V(D_{i})\rightarrow C[x_{a+1},x_{b-1}])=\emptyset.
$$
So,  there are vertices $x_a,x_b\in C_m$ ($x_a\not= x_b$), $u\in V(D_k)$ and $v\in V(D_p)$ ($1\leq k \leq p \leq s$) such that 
$$
x_au,\, vx_b\in A(D) \quad \hbox {and} \quad A(R \rightarrow \cup_{i=1}^p V(D_{i})) = A(\cup_{i=k}^s V(D_i) \rightarrow R)=\emptyset, \eqno (8)
$$
where $R:=V(C_m[x_a,x_b]) \setminus \{x_a,x_b\}$. In particular, $A(R,\cup_{i=k}^p V(D_i))=\emptyset$. It is clear that $|R|\not=0$. We extend the path $C_m[x_b,x_a]$ with the vertices of $R$ as much as possible. 
Let $Q$ be a obtained extended path. From the maximality of $C_m$ it follows that some vertices $y_1,y_2, \ldots , y_d\in R$ ($1\leq d \leq |R|$) do not on the extended path $Q$. Let $y_i\in \{y_1,y_2,\ldots ,y_d\}$ and $z\in V(D_t)$ with $t\in [k,p]$ are arbitrary vertices. Using Lemma 4.4 and (8) we obtain, 
$$
d(y_i,V(C_m))\leq m+d-1; \quad d(z,V(C_m))\leq m-|R|+1; \eqno (9)
$$
$$
 d(y_i,B)\leq \left\{ \begin{array}{lc} n-m-d_t,\quad \hbox{if}\quad \langle B \rangle \quad \hbox{is unilaterally connected},    \\ n-m-d_t-1+d(y_i,\{z_0\}),  \quad \hbox{otherwise } ; \\ \end{array} \right. \eqno (10)
$$
$$
 d(z,B)\leq \left\{ \begin{array}{lc} n-m+d_t-2,\quad \hbox{if}\quad \langle B \rangle \quad \hbox{is unilaterally connected},    \\ n-m+d_t-3+d(z,\{z_0\}),  \quad \hbox{otherwise} ; \\ \end{array} \right. \eqno (11)
$$
where $d_i:=|V(D_i)|$.

Suppose first that $\langle B \rangle$ is not unilaterally connected. Then by Remark 4.6, $z_0\in B$ and $d(w,\{z_0\})\leq 1$ for all $w\in B$. From (9)-(11) it follows that
$$
2n-1\leq d(y_i)+d(z) \leq m+d-1+n-m-d_t-1+d(y_i,\{z_0\})+m-|R|+1 +$$
 $$n-m+d_t-3+d(z,\{z_0\})\leq  
 2n+d-|R|-4+d(z,\{z_0\})+d(y_i,\{z_o\}).   
$$
Since $d\leq |R|$, $d(z,\{z_0\})\leq 1$ and $d(y_i,\{z_0\})\leq 2$ we see that  $d= |R|$, $d(z,\{z_0\})= 1$, $d(y_i,\{z_0\})=2$ and $d(z,V(C_m))=m-|R|+1$. From $d(y_i,\{z_0\})=2$ ($y_i$ is arbitrary) and the maximality of $C_m$ it follows that $d=|R|=1$ ($R=\{x_{a+1}\}$), $u=v$, i.e., $k=p$ and $ux_{a+2}\in A(D)$. Since $z$ is an arbitrary vertex, we have that $uz_0\in A(D)$ or $z_0u\in A(D)$ (when $z=u$). Therefore, if $uz_0\in A(D)$, then $C_{m+2}=x_auz_0x_{a+1}\ldots x_a$, and if $z_0u\in A(D)$, then $C_{m+2}=x_ax_{a+1}z_0ux_{a+2}\ldots x_a$, which contradicts the assumption that $C_m$ is a cycle of the maximum length.

Suppose now that $\langle B \rangle$ is  unilaterally connected. By (9)-(11) we have
$$
d(y_i)+d(z) \leq m+d+m-|R|+n-m-d_t+n-m+d_t-2 = 2n+d-|R|-2\leq 2n-2. \eqno {(12)} 
$$

First consider the case $z_0\in B$. From (12) it follows that $V(D_t)=\{z_0\}$ ($k=p$) and $s\geq 2$ since $n-m\geq 2$. 

Let $z_0\notin V(D_s)$, i.e., $t<s$. If $A(R\rightarrow V(D_s))=\emptyset$, then by (8) we have $A(R,V(D_s))=\emptyset$. Then similarly as in (12), we can  show that  for each $x\in V(D_s)$  and $y_i $, $i\in [1,d]$ the following holds
$$
2n-1\leq d(x)+d(y_i)\leq 2n+d-|R|-3\leq 2n-3$$
 since $d\leq |R|$ and $d_s\geq 1$, which is a contradiction.
 So, we can assume that $A(R\rightarrow V(D_s))\not=\emptyset$. Recall that 
 $$
A(V(D_s)\rightarrow V(C_m))\not= \emptyset  \quad \hbox{and} \quad A(V(D_s)\rightarrow R)=\emptyset
$$ 
(by (7) and (8), respectively). 
Therefore, there are two vertices $v_1\in V(D_s)$ and $x_q\in V(C_m[x_b,x_a])$ such that $A(V(D_s)\rightarrow V(C_m[x_{a+1},x_{q-1}]))= \emptyset$ and $v_1x_q\in A(D)$.
 Now we can claim that there are two vertices $x_r\in V(C_m[x_{a+1},x_{q-2}])$ and $u_1\in B$  such that $x_ru_1\in A(D)$ and  $A(V(C_m[x_{r+1},x_{q-1}])\rightarrow B)= \emptyset$. In particular, we have $A(V(D_s),F)=\emptyset$, where $F:=V(C_m[x_{r+1},x_{q-1}])$.
 We extend the path $C_m[x_{q},x_{r}]$ with vertices $x_{r+1},\ldots , x_{q-1}$ as much as possible. Then some vertices $u_1,u_2,\ldots , u_j$ of $F$ ($1\leq j\leq |F|$) do not on the obtained extended path. Therefore, by Lemma 4.4, for all vertices $u_i\in F$ and $x\in V(D_s)$ the following hold
$$d(x,V(C_m))\leq m-|F|+1;\quad d(x,B)\leq n-m+d_s-2;$$
 $$d(u_i,V(C_m))\leq m+j-1; \quad d(u_i,B)\leq n-m-d_s.$$ 
Hence 
$d(x)+d(u_i)\leq 2n-2$,  a contradiction since $z_0\notin \{u_i,x\}$ and the vertices $u_i$, $x$ are nonadjacent.

 Let now $z_0\in V(D_s)$. Then $t=s$ and  $V(D_s)=\{z_0\}$ ($s\geq 2$). It is not difficult to see that for the   converse  digraph of $D$ we have the considered case $z_0\notin V(D_s)$.

Now consider the case $z_0\notin B$, i.e., $z_0\in V(C_m)$. Then from (12) it follows that $d=1$, $y_1=z_0$ and $u=v$ (i.e., there is a path, say $P$, from $x_b$ to $x_a$ with vertex set $V(C_m)\setminus \{z_0\}$. Therefore, $uPu$ is a cycle of length $m$, which does not contain the vertex $z_0$. So, we have one of  above-considered cases. This completes the proof of Theorem 4.5 . \fbox \\\\

Using Theorem 4.5 and Lemma 4.4 it is not difficult to show that the following corollaries are true.\\

\noindent\textbf{Corollary 4.7}. {\it Let $D$ be a strong digraph of order $n\geq 3$ satisfying the conditions of Theorem 4.5. Then $D$ has a cycle that contains all the vertices of $D$ maybe except $z_0$.}\\

\noindent\textbf{Corollary 4.8}. {\it Let $D$ be a strong digraph of order $n\geq 3$. If in $D$ the degrees of $n-1$ vertices  at least $n$, then $D$ is a Hamiltonian or contains a cycle of length $n-1$ (in fact $D$ has a cycle that contains all the vertices with degree at least $n$).} \\

 \noindent\textbf{Note added in proof} (for section 4). 
Later on, Berman and Liu \cite{[16]} and  Li, Flandrin and Shu \cite{[17]} proved Theorems 4.9 and 4.10, respectively, which improved Theorem 4.5. Before to formulate these theorems we need the following definitions. 

Let $D$ be a digraph of order $n\geq 3$ and let $M$ be a non-empty subset of $V(D)$. Following \cite {[16]} and \cite {[17]}, we say that

(i) the subset $M$ is Meyniel set if $d(x)+d(y)\geq 2n-1$ for every pair of vertices $x, y$ of $M$ which are non-adjacent in
$D$.

(ii) the digraph $D$ is $M$-strongly connected if for any pair of distinct vertices $x,y$ of $M$ there exists a path from $x$ to $y$ and a path from $y$ to $x$ in $D$.\\

\noindent\textbf{Theorem 4.9} (Berman and Liu \cite{[16]}). {\it Let $D$ be a strong digraph. If   $M$ is a Meyniel set in $D$, then $D$ has a cycle containing  all the vertices of $M$.}
\\

\noindent\textbf{Theorem 4.10} (Li, Flandrin and Shu \cite{[17]}). {\it Let $D$ be a digraph of order $n$ and $M$ be a Meyniel set in $D$. If $D$ is $M$-strongly connected, then $D$ has a cycle through all the vertices of $M$.}

\end{document}